\documentclass[11pt,reqno]{amsart}
\usepackage[T1]{fontenc}
\usepackage{lmodern}
\usepackage{amsmath,amssymb,amsthm}
\usepackage[margin=1.1in]{geometry}
\usepackage{microtype}
\usepackage{booktabs}
\usepackage{hyperref}
\hypersetup{colorlinks=true,linkcolor=blue,citecolor=blue,urlcolor=blue,
 pdftitle={Spine extremality and its frequency in critical branching random walks},
 pdfauthor={Haoming Zhang}}
\newcommand{\Pp}{\mathbf P}
\newcommand{\Q}{\mathbf Q}
\newcommand{\E}{\mathbf E}
\newcommand{\EQ}{\mathbf E_{\mathbf Q}}
\newcommand{\1}{\mathbf 1}

\newcommand{\R}{\mathbb R}
\newcommand{\cT}{\mathcal T}
\newcommand{\cF}{\mathcal F}
\newcommand{\cH}{\mathcal H}

\newtheorem{theorem}{Theorem}[section]
\newtheorem{proposition}[theorem]{Proposition}
\newtheorem{lemma}[theorem]{Lemma}
\newtheorem{corollary}[theorem]{Corollary}
\theoremstyle{remark}
\newtheorem{remark}[theorem]{Remark}
\numberwithin{equation}{section}

\title[Spine extremality and its frequency]
{Spine extremality and its frequency in critical branching random walks}
\author{Haoming Zhang}
\date{16 October 2026}

\begin{document}
\begin{abstract}
We study the generations at which the distinguished spine of a critical
branching random walk occupies the spatial minimum. Under finite positive
offspring variance and any nonconstant law of independent real edge
increments, their number up to generation $N$ is
$(\log N)^{1+o(1)}$ almost surely and has order $\log N$ in probability,
with ties included. The same conclusions hold for the unique-minimum
event. No spatial moment assumption is needed. The proof combines
counting estimates from the spine decomposition with uniform control of
extremal multiplicity in an ordinary surviving tree. For a selected
minimizer, its one-time probability is the ordinary survival probability
$q_n$, and recurrence is equivalent to $\sum_nq_n=\infty$.
We also exhibit a critical infinite-variance offspring family for which
the weak-minimum event occurs only finitely often, for every nonconstant
increment law. We determine the polynomial order of its last occurrence
and, for Gaussian or Laplace increments, prove quantitative eventual
separation of the spine from the minimum. The sharp frequency constant
and limiting distribution in the nonconstant finite-variance model
remain undetermined.
\end{abstract}
\maketitle

\section{Introduction and main results}

A critical, nondegenerate Galton--Watson tree is finite almost surely.
Its size-biased version has a distinguished infinite ancestral line, the
\emph{spine}, with ordinary critical trees attached to it. Assigning
independent displacements to the edges gives a branching random walk
which survives at every generation. The question proposed by Antal
J\'arai asks whether the spine particle is at the spatial minimum in
finitely or infinitely many generations, and, in the latter case, how
frequently this occurs. We prove recurrence and determine the almost-sure
logarithmic exponent of the number of such generations in the
finite-variance setting, including spatial ties. We also study a family
of infinite-variance offspring laws for which only finitely many such
generations occur, and obtain spatial separation estimates in that regime.

Kesten's incipient infinite tree construction~\cite{K86} provides the
genealogical setting. The oriented version embeds a vertex $u$ at
$(V(u),|u|)$ in space--time. J\'arai and Nachmias~\cite{JN14} study
electrical resistance of this infinite oriented branching random walk.
Here particles retain their genealogical identities even when their
space--time locations coincide.

There are two distinct notions of a leftmost particle. An ordered family
tree has a genealogically leftmost vertex at each nonempty level, whereas
a branching random walk has one or more particles of smallest spatial
position. These need not coincide. The size-biased construction of
Lyons, Pemantle and Peres~\cite{LPP95} provides the basic change of measure
used here. Their proof of the critical survival estimate uses the
genealogically leftmost spine event. A related conditioning argument
appears in Engl\"ander and Peres~\cite{EP17}, again with genealogical
leftmost order. The present question concerns spatial order over
infinitely many generations.

Work on extremal displacement, including Kesten~\cite{K95} and
Lalley and Shao~\cite{LS15}, concerns the position of the population's
extreme rather than the generations at which a distinguished immortal
line occupies that extreme. This distinction also separates the present
counting problem from a fixed-time conditional limit theorem.

Our first observation is that a spine particle which is a minimum of
the whole generation is also a minimum among the descendants of any
earlier spine particle. Independence along the spine then gives
multipoint upper bounds and a stopping-time lower-tail estimate for the
count. To include ties, we bound the conditional second moment of the
number of minimizers. The counting estimates apply to both weak and strict minima once their
one-time probabilities are comparable to the ordinary survival probability.
For weak minima this comparison follows from the multiplicity bound; for
strict minima we compare the corresponding recursive measures directly.

\subsection{The model}

Let $p=(p_k)_{k\geq0}$ be an offspring distribution satisfying
\begin{equation}\label{eq:critical}
 \sum_{k\geq0}kp_k=1,\qquad p_1<1.
\end{equation}
Its probability generating function is
$f(s)=\sum_{k\geq0}p_ks^k$, $0\leq s\leq1$.
Write $\Pp$ for the law of an ordered Galton--Watson tree $\cT$ with this
offspring distribution, started from a single root. Vertices are labelled
by finite words in positive integers: the root is $\varnothing$, and the
children of $u$ are $uj$, $1\leq j\leq k(u)$, where $k(u)$ is its offspring
count. Put
$\cT_n=\{u\in\cT:|u|=n\}$ and $Z_n=|\cT_n|$.
Each edge carries a real displacement with law $\mu$, independently of
the tree and all other displacements. The position $V(u)$ is the sum of
the displacements along the path from the root to $u$, with
$V(\varnothing)=0$.

The marked size-biased law $\Q$ is defined as follows. The root $w_0$ is
on the spine. Each spine vertex $w_i$ has a number $\widehat K_i$ of
children with distribution
\begin{equation}\label{eq:sizebias}
 \Q(\widehat K_i=k)=kp_k,\qquad k\geq1.
\end{equation}
One child, chosen uniformly, is $w_{i+1}$. Every other child starts an
independent ordinary Galton--Watson tree with offspring law $p$.
All choices at successive spine vertices are independent, and every
edge still has an independent displacement of law $\mu$.
We use $\Q$ for the joint law of the tree, its spatial marks and its spine.
In particular, $Z_n\geq1$ under $\Q$.

For either law, write
\[
 M_n=\min_{u\in\cT_n}V(u),\qquad
 K_n=\#\{u\in\cT_n:V(u)=M_n\},
\]
with $M_n=+\infty$ and $K_n=0$ when $Z_n=0$. Thus $M_n$ records the
minimum position, whereas $K_n$ counts the particles attaining it.
For each nonempty generation, let $\ell_n$ be the lexicographically
first of these minimizers. Selecting one minimizer gives an event whose
probability is the ordinary survival probability, even when several
particles tie. Define
\begin{align}
 B_n&=\{w_n=\ell_n\},\label{eq:selected-event}\\
 A_n&=\{V(w_n)=M_n\}.\label{eq:weak-event}
\end{align}
Thus $B_n\subseteq A_n$; the second event includes ties. Write
\begin{equation}\label{eq:q-and-count}
 q_n=\Pp(Z_n>0),\qquad
 C_N=\sum_{n=1}^N\1_{B_n},\qquad
 s_N=\sum_{n=1}^Nq_n.
\end{equation}
For the event including ties, put
\begin{equation}\label{eq:weak-count}
 a_n=\Q(A_n),\qquad
 D_N=\sum_{n=1}^N\1_{A_n},\qquad t_N^A=\sum_{n=1}^Na_n.
\end{equation}
We also write $D_\infty=\lim_{N\to\infty}D_N$, allowing the value
$+\infty$.
Expectations under $\Pp$ and $\Q$ are denoted by $\E$ and $\EQ$,
respectively. All logarithms are natural.
For positive sequences, $b_n\asymp c_n$ means that $b_n/c_n$ is
bounded above and below by positive constants for all sufficiently
large $n$.

\begin{center}
\small
\begin{tabular}{@{}ll@{}}
\toprule
Notation & Meaning \\
\midrule
$\Pp,\ \Q$ & Ordinary law and marked size-biased law \\
$Z_n,\ q_n$ & Population size and ordinary survival probability \\
$V(u),\ w_n$ & Particle position and spine vertex \\
$M_n,\ K_n$ & Spatial minimum and number of minimizers \\
$A_n,\ B_n$ & Spine is a minimizer; spine is the selected minimizer \\
$D_N,\ C_N$ & Counts of $A_n$ and $B_n$ up to generation $N$ \\
$t_N^A,\ s_N$ & Expected counts under $\Q$ \\
\bottomrule
\end{tabular}
\end{center}

\subsection{Main results}

The principal finite-variance conclusion answers the spatial question,
including ties. The offspring-variance assumption is
\begin{equation}\label{eq:variance}
 0<\sigma^2=\sum_{k\geq0}(k-1)^2p_k<\infty.
\end{equation}
The displacement law need not be centered or have any finite moment.

\begin{theorem}[Frequency with and without ties]\label{thm:frequency-main}
Assume \eqref{eq:critical} and \eqref{eq:variance}. For every $\mu$,
\begin{equation}\label{eq:selected-exponent}
 \frac{\log C_N}{\log\log N}\longrightarrow1
 \qquad\Q\text{-almost surely}.
\end{equation}
If $\mu$ is not a point mass, there is a finite constant $K$ such that
\begin{equation}\label{eq:weak-comparison}
 q_n\leq a_n\leq Kq_n\quad(n\geq1),\qquad
 \frac{\log D_N}{\log\log N}\longrightarrow1
 \quad\Q\text{-almost surely}.
\end{equation}
For nonconstant $\mu$, in particular,
$D_N=(\log N)^{1+o(1)}$ almost surely and $D_N/N\to0$.
Moreover, $C_N$ has order $\log N$ in probability for every $\mu$;
the same holds for $D_N$ when $\mu$ is nonconstant. Here a sequence
$W_N$ has \emph{order $b_N$ in probability} if, for every
$\varepsilon>0$, there are constants
$0<c_\varepsilon<C_\varepsilon<\infty$ such that
\[
 \liminf_{N\to\infty}
 \Q(c_\varepsilon b_N\leq W_N\leq C_\varepsilon b_N)
 \geq1-\varepsilon.
\]
For nonconstant $\mu$, put $\tau_k=\inf\{N:D_N\geq k\}$.
Then all $\tau_k$ are finite and
\begin{equation}\label{eq:arrival-exponent}
 \frac{\log\log\tau_k}{\log k}\longrightarrow1
 \qquad\Q\text{-almost surely}.
\end{equation}
If instead $\mu=\delta_c$, then $A_n$ holds for every $n$ and $D_N=N$.
\end{theorem}

Theorem~\ref{thm:frequency-main} is proved after the general counting
estimates and the treatment of ties. Section~\ref{sec:strict-minimum}
proves the same recurrence and frequency exponent for nonconstant $\mu$
when the spine is required to be the \emph{unique} minimum.
The logarithmic exponent and the order in probability do not assert
that $D_N/\log N$ converges.

The selected event $B_n$ also admits a criterion without any offspring
variance assumption. Its role is to connect extremality to the ordinary
survival probabilities.

\begin{theorem}[A recurrence criterion]\label{thm:criterion}
Assume \eqref{eq:critical}, and let $\mu$ be any probability law on $\R$.
For every $n\geq1$ and $1\leq m<n$,
\begin{equation}\label{eq:identities-main}
 \Q(B_n)=q_n,
 \qquad \Q(B_m\cap B_n)\leq q_mq_{n-m}.
\end{equation}
Moreover,
\begin{equation}\label{eq:criterion}
 \Q(B_n\text{ infinitely often})=
 \begin{cases}
 1,&\displaystyle\sum_{n=1}^{\infty}q_n=\infty,\\[3pt]
 0,&\displaystyle\sum_{n=1}^{\infty}q_n<\infty.
 \end{cases}
\end{equation}
For each $N\geq1$,
\begin{equation}\label{eq:moments}
 \EQ[C_N]=s_N,
 \qquad \EQ[C_N^2]\leq s_N+2s_N^2.
\end{equation}
\end{theorem}

Recurrence in the finite-variance case is an immediate consequence of the criterion
and the classical critical survival estimate.

\begin{corollary}[Finite-variance spatial extremality]\label{cor:main}
Assume \eqref{eq:critical} and \eqref{eq:variance}.
Then, for every displacement law $\mu$,
\begin{equation}\label{eq:io-main}
 \Q(B_n\text{ infinitely often})
 =\Q(A_n\text{ infinitely often})=1.
\end{equation}
Furthermore,
\begin{equation}\label{eq:asymptotics}
 \Q(B_n)\sim\frac{2}{\sigma^2 n},\qquad
 \EQ[C_N]\sim\frac{2}{\sigma^2}\log N.
\end{equation}
If $\mu$ is atomless, then $A_n=B_n$ almost surely for every $n$, so
\eqref{eq:asymptotics} also holds for the events $A_n$ and their count.
\end{corollary}

\begin{remark}
The variance $\sigma^2$ in \eqref{eq:asymptotics} is the offspring
variance. No centering, variance or other moment assumption is imposed
on $\mu$. Replacing every edge displacement by its negative shows that
the spine is also at the spatial maximum infinitely often almost surely. This does
not assert that the minimum and maximum events occur simultaneously.
\end{remark}

\begin{theorem}[An infinite-variance finite-occurrence regime]
\label{thm:finite-main}
Suppose that the offspring generating function is
\[
 f(s)=s+c(1-s)^\alpha,\qquad
 1<\alpha<2,\qquad 0<c\leq1/\alpha,
\]
and $\mu$ is nonconstant. Put $\beta=1/(\alpha-1)>1$, and let
$\mathcal L=\sup\{n\geq1:A_n\text{ occurs}\}$, with
$\sup\varnothing=0$. Then
\[
 q_n\sim[c(\alpha-1)n]^{-\beta},\qquad a_n\asymp n^{-\beta},
 \qquad \EQ[D_\infty]<\infty.
\]
In particular, $\mathcal L<\infty$ almost surely, and
\[
 \Q(\mathcal L>N)\asymp N^{1-\beta},\qquad
 \EQ[\mathcal L^u]<\infty\quad\Longleftrightarrow\quad u<\beta-1
 \quad(u>0).
\]
For Gaussian or Laplace increments, the distance from the spine to the
minimum tends to infinity almost surely, with the quantitative bounds
of Theorem~\ref{thm:quantitative-gap}.
\end{theorem}

The generating-function survival asymptotic is a special case of the
classical infinite-variance theory of Slack~\cite{Slack68}; we give its
short direct proof for the displayed family. The comparison of $a_n$
with $q_n$ requires a separate treatment of spatial ties.
Theorem~\ref{thm:finite-main} follows from
Proposition~\ref{prop:stable-multiplicity},
Proposition~\ref{prop:last-occurrence-tail}, and
Theorem~\ref{thm:quantitative-gap} below.

These results distinguish the two alternatives in the original question.
Under finite positive offspring variance, the finite-occurrence branch
does not arise. Without that variance assumption, either behavior can
occur. We do not give a classification of every infinite-variance
critical offspring law, nor an exact spatial separation rate in the
finite-occurrence regime.

\subsection{Organization of the proof}

Sections~\ref{sec:finite-identities}--\ref{sec:recurrence} establish
the size-biased identities, the independent spine decomposition, and the
survival-series criterion. Section~\ref{sec:generic-frequency} develops
counting estimates in terms of the expected number of occurrences.
Sections~\ref{sec:multiplicity-bound} and~\ref{sec:strict-minimum}
supply the spatial input by controlling ties and unique minima.
Sections~\ref{sec:consequences} and~\ref{sec:fluctuations} discuss arrival
times and fluctuations. Section~\ref{sec:stable-multiplicity} treats
finite occurrence and eventual spatial separation, and
Section~\ref{sec:questions} records the remaining questions.

\section{Finite-generation identities}
\label{sec:finite-identities}

Let $\cF_n$ be the sigma-field generated by the tree through generation
$n$ and its edge displacements through that generation, without recording
the distinguished spine. Reproduction of vertices in generation $n$
is not included in $\cF_n$.

\begin{lemma}[The finite-generation spine identity]\label{lem:spine}
For every nonnegative measurable function $H$ of a spatial tree through
generation $n$ and one of its generation-$n$ vertices,
\begin{equation}\label{eq:many-one}
 \EQ\bigl[H(\cT_{\leq n},V,w_n)\bigr]
 =\E\left[\sum_{u\in\cT_n}H(\cT_{\leq n},V,u)\right].
\end{equation}
In particular,
\begin{equation}\label{eq:density-uniform}
 \left.\frac{d\Q}{d\Pp}\right|_{\cF_n}=Z_n,
 \qquad
 \Q(w_n=u\mid\cF_n)=\frac{1}{Z_n}
 \quad (u\in\cT_n).
\end{equation}
\end{lemma}

\begin{proof}
Fix a finite ordered tree through generation $n$ and a vertex $u$ at
level $n$. At every ancestor of $u$ having $k$ children, the factor
contributed under $\Q$ by its offspring count and its specified spine
child is $(kp_k)/k=p_k$. At every other vertex below level $n$,
the offspring factor is also $p_k$. Hence the probability of the fixed
tree together with the marked path to $u$ is its ordinary tree
probability under $\Pp$. Conditional edge-displacement laws are the
same under both measures. Integrating over the displacements and
summing over all such trees and vertices proves \eqref{eq:many-one},
first for cylinder functions and then for nonnegative measurable $H$.
Summing over the marked vertex yields the density in
\eqref{eq:density-uniform}; equality of the contributions of the
generation-$n$ vertices gives conditional uniformity.
\end{proof}

\begin{remark}
The conditioning in \eqref{eq:density-uniform} is on the finite tree
through generation $n$. It is not conditioning on the entire infinite
tree, whose unique infinite ancestral line identifies the spine.
The underlying size-bias identity is classical; see~\cite[Section 2]{LPP95}.
\end{remark}

\begin{proposition}[Conditioning on distant survival]
\label{prop:conditioning-limit}
For every fixed $m$, the law of the unmarked spatial tree through
generation $m$, under $\Pp(\,\cdot\mid Z_N>0)$, converges in total
variation to its law under $\Q$ as $N\to\infty$.
\end{proposition}
\begin{proof}
By the branching property, its density relative to $\Pp$ on $\cF_m$ is
\[
 R_{m,N}=\frac{1-(1-q_{N-m})^{Z_m}}{q_N}.
\]
Criticality and nondegeneracy give $q_n\downarrow0$. Since $f'(1)=1$,
the recursion $q_{n+1}=1-f(1-q_n)$ gives $q_{n+1}/q_n\to1$.
Consequently $q_{N-m}/q_N\to1$ for fixed $m$, and
$R_{m,N}\to Z_m$ almost surely. Also
$0\leq R_{m,N}\leq(q_{N-m}/q_N)Z_m\leq2Z_m$ for all sufficiently
large $N$. As $\E[Z_m]=1$, dominated convergence gives convergence
of the densities in $L^1(\Pp)$. The density $Z_m$ is exactly the
unmarked $\Q$-density in \eqref{eq:density-uniform}.
\end{proof}

\begin{proposition}[Selection and multiplicity]\label{prop:selection}
For each $n\geq1$,
\begin{equation}\label{eq:selected-prob}
 \Q(B_n)=q_n.
\end{equation}
Recall that $K_n$ is the number of spatial minimizers, with $K_n=0$ on extinction. Then
\begin{equation}\label{eq:multiplicity}
 \Q(A_n)=\E[K_n].
\end{equation}
\end{proposition}

\begin{proof}
For a nonempty generation there is exactly one selected vertex $\ell_n$.
Applying \eqref{eq:many-one} with $H(t,x,u)=\1_{\{u=\ell_n(t,x)\}}$
gives \eqref{eq:selected-prob}. Applying it instead with $H$ equal to
the indicator that $u$ is a spatial minimizer gives
\eqref{eq:multiplicity}.
Equivalently, conditional uniformity gives
$\Q(B_n\mid\cF_n)=1/Z_n$; multiplying by the density $Z_n$
exactly cancels size-biasing.
\end{proof}

\begin{lemma}[Absence of spatial ties]\label{lem:ties}
If $\mu$ is atomless, distinct vertices of any fixed generation have
distinct positions almost surely, under both $\Pp$ and $\Q$.
\end{lemma}

\begin{proof}
Condition on the finite genealogical tree, and consider distinct
vertices $u,w$ at the same level. Their paths separate after their
last common ancestor. Choose an edge on the path to $u$ after this
separation, and condition on all other edge displacements. Then
$V(u)-V(w)$ is that edge displacement plus a fixed constant.
It is zero with probability zero. There are finitely many pairs
in each generation, and countably many generations.
\end{proof}

\section{Dependence between generations and the tail event}

\subsection{Independent data along the spine}

The size-biased construction can be encoded by independent, identically
distributed blocks $(\mathcal D_i)_{i\geq0}$. A block $\mathcal D_i$
records the offspring count of $w_i$, the index of its spine child,
the displacements to all its children, and the complete ordinary spatial
trees attached to its nonspine children. The roots of these attached
trees are recorded relative to their attachment vertex. Their labels
are also relative labels. All blocks have the same law.

Define
\begin{equation}\label{eq:past-blocks}
 \cH_m=\sigma(\mathcal D_0,\ldots,\mathcal D_{m-1}).
\end{equation}
This sigma-field records more than the population through time $m$,
because it includes the future of earlier off-spine trees. Nevertheless,
the whole marked descendant tree rooted at $w_m$, after translating
$V(w_m)$ to zero and deleting the common label prefix $w_m$, depends
only on $\mathcal D_m,\mathcal D_{m+1},\ldots$. It therefore has law
$\Q$ and is independent of $\cH_m$. The event $B_m$ belongs to $\cH_m$.

For $k\geq1$, let $B_k^{(m)}$ denote the selected-minimum event at
relative generation $k$ in this descendant tree.

\begin{lemma}[Two-time estimate]\label{lem:twotime}
For $1\leq m<n$,
\begin{equation}\label{eq:two-bound}
 \Q(B_m\cap B_n)\leq q_mq_{n-m}.
\end{equation}
\end{lemma}

\begin{proof}
If $B_n$ occurs, then $w_n$ minimizes position among all generation-$n$
vertices and hence among those descended from $w_m$. If there were a
lexicographically earlier spatial minimizer in the latter subtree,
it would also contradict $B_n$ in the whole tree. Removing the shared
label prefix preserves lexicographic order, and subtracting $V(w_m)$
preserves spatial order. Consequently,
\begin{equation}\label{eq:event-inclusion}
 B_n\subseteq B_{n-m}^{(m)}.
\end{equation}
Since $B_m\in\cH_m$ and $B_{n-m}^{(m)}$ is independent of $\cH_m$,
Proposition~\ref{prop:selection} gives
\[
 \Q(B_m\cap B_n)
 \leq\Q(B_m\cap B_{n-m}^{(m)})
 =\Q(B_m)\Q(B_{n-m}^{(m)})=q_mq_{n-m}.
\]
\end{proof}

\begin{lemma}[Zero--one law]\label{lem:zeroone}
Under \eqref{eq:critical},
\begin{equation}\label{eq:zeroone}
 \Q(B_n\text{ infinitely often})\in\{0,1\}.
\end{equation}
\end{lemma}

\begin{proof}
Let $E=\{B_n\text{ infinitely often}\}$ and fix $m\geq1$.
Only finitely many off-spine trees attach to
$w_0,\ldots,w_{m-1}$, since every offspring count is finite.
Each attached ordinary tree becomes extinct almost surely under
\eqref{eq:critical}. Thus there is an almost surely finite random
generation $R_m\geq m$ beyond which none of these trees has any vertices.
For every $n>R_m$, all particles in generation $n$ descend from $w_m$.
Translation of all positions by a common constant, and deletion of
a common label prefix, preserve the selected minimum. Hence
\[
 \1_{B_n}=\1_{B_{n-m}^{(m)}}\qquad(n>R_m)
\]
almost surely. It follows that $E$ agrees almost surely with
\[
 E_m=\{B_k^{(m)}\text{ infinitely often as }k\to\infty\},
\]
which belongs to $\sigma(\mathcal D_i:i\geq m)$ and is independent
of $\cH_m$.

Therefore $E$ is independent of $\cH_m$ for every $m$. The increasing
union of these sigma-fields generates the sigma-field of the full
sequence of blocks, which also contains $E$. A monotone-class argument
shows that $E$ is independent of that entire sigma-field, and therefore
of itself. Thus $\Q(E)=\Q(E)^2$, proving the assertion.
\end{proof}

\section{The recurrence criterion and finite-variance survival}
\label{sec:recurrence}

\begin{proof}[Proof of Theorem~\ref{thm:criterion}]
The identities and inequalities in \eqref{eq:identities-main} have been
proved in Proposition~\ref{prop:selection} and Lemma~\ref{lem:twotime}.
Summing the one-time identity gives $\EQ[C_N]=s_N$. Expanding the square
and applying the two-time estimate gives
\begin{align}
 \EQ[C_N^2]
 &=s_N+2\sum_{1\leq m<n\leq N}\Q(B_m\cap B_n)\notag\\
 &\leq s_N+2\sum_{m=1}^{N-1}\sum_{k=1}^{N-m}q_mq_k
 \leq s_N+2s_N^2.\label{eq:second-moment}
\end{align}
This proves \eqref{eq:moments}.

If $\sum_n q_n<\infty$, the first Borel--Cantelli lemma implies
that only finitely many $B_n$ occur almost surely.
Suppose instead that $s_N\to\infty$. For fixed $r\geq1$, put
$C_{r,N}=\sum_{n=r}^N\1_{B_n}$, and use $s_0=0$.
The Cauchy--Schwarz inequality yields
\begin{align}
 \Q\left(\bigcup_{n=r}^N B_n\right)
 &=\Q(C_{r,N}>0)
 \geq\frac{(\EQ[C_{r,N}])^2}{\EQ[C_{r,N}^2]}\notag\\
 &\geq\frac{(s_N-s_{r-1})^2}{s_N+2s_N^2}.
 \label{eq:positive}
\end{align}
Letting $N\to\infty$ and then $r\to\infty$ gives
$\Q(B_n\text{ infinitely often})\geq1/2$.
This is the second-moment implication underlying the
Kochen--Stone bound~\cite{KS64}; \eqref{eq:positive} supplies the
argument directly in the present setting. Lemma~\ref{lem:zeroone}
upgrades positive probability to probability one.
\end{proof}

For completeness, the survival asymptotic needed for the corollary
is recalled with a short generating-function proof. It is the classical
Kolmogorov estimate under the finite-variance hypothesis; see
Kesten, Ney and Spitzer~\cite{KNS66} and~\cite[Theorem C]{LPP95}.

\begin{lemma}[Critical survival]\label{lem:survival}
Under \eqref{eq:critical} and \eqref{eq:variance},
\begin{equation}\label{eq:survival}
 q_n\sim\frac{2}{\sigma^2 n}.
\end{equation}
\end{lemma}

\begin{proof}
Let $f(s)=\sum_{k\geq0}p_ks^k$ be the offspring generating function.
Criticality and nondegeneracy imply that the extinction probability is
one: indeed $f(s)>s$ for $0\leq s<1$, so the only fixed point in
$[0,1]$ is one. Thus $q_n\downarrow0$.
Also $f'(1)=1$ and $f''(1)=\sigma^2$. Taylor expansion from the left gives
\[
 f(1-x)=1-x+\frac{\sigma^2}{2}x^2+o(x^2)
 \qquad(x\downarrow0).
\]
The branching property implies
\[
 q_{n+1}=1-f(1-q_n)
 =q_n-\frac{\sigma^2}{2}q_n^2+o(q_n^2).
\]
In particular, $q_{n+1}/q_n\to1$, and consequently
\[
 \frac{1}{q_{n+1}}-\frac{1}{q_n}
 =\frac{q_n-q_{n+1}}{q_nq_{n+1}}
 \longrightarrow\frac{\sigma^2}{2}.
\]
Taking arithmetic means of these increments yields
$1/(nq_n)\to\sigma^2/2$, proving \eqref{eq:survival}.
\end{proof}

\begin{proof}[Proof of Corollary~\ref{cor:main}]
Lemma~\ref{lem:survival} implies $\sum_nq_n=\infty$.
Theorem~\ref{thm:criterion} therefore gives infinitely many occurrences
of $B_n$ almost surely, and $B_n\subseteq A_n$ gives the same conclusion
for $A_n$. The probability asymptotic follows from $\Q(B_n)=q_n$.
Summing \eqref{eq:survival} gives
$s_N\sim(2/\sigma^2)\log N$, which proves the expectation asymptotic.
When $\mu$ is atomless, Lemma~\ref{lem:ties} identifies $A_n$ and $B_n$
simultaneously for all $n$ outside a null set.
\end{proof}

\section{Occupation counts in their expected time scale}
\label{sec:generic-frequency}

The same argument controls both spatial-minimum counts and
selected-minimum counts. We formulate it once in terms of their expected
number of occurrences. Throughout this section, choose either
$E_n=A_n$ for all $n$ or $E_n=B_n$ for all $n$, and set
\begin{equation}\label{eq:generic-count-notation}
 r_n=\Q(E_n),\qquad
 U_N=\sum_{n=1}^N\1_{E_n},\qquad
 t_N=\sum_{n=1}^Nr_n=\EQ[U_N],\qquad t_0=0.
\end{equation}
No assumption on ties is needed: each choice of events has its own
sequence $(r_n)$ and expected count $(t_N)$.

The argument has three parts. Restriction to a descendant tree gives
bounds on joint probabilities and hence exponential moments. Applying
the same restriction at the time of the $k$th occurrence yields a
lower-tail estimate. We then apply Borel--Cantelli on a geometric grid
in $t_N$ to obtain almost-sure bounds for the counts.

Write $E_k^{(m)}$ for the corresponding event at relative generation $k$
in the descendant tree rooted at $w_m$. The independent spine blocks give
\begin{equation}\label{eq:generic-hereditary}
 E_n\in\cH_n,\qquad E_n\subseteq E_{n-m}^{(m)}\quad(0\leq m<n),
\end{equation}
and $E_k^{(m)}$ is independent of $\cH_m$, with probability $r_k$.
For $A_n$, the inclusion holds because a spatial minimum remains a
minimum after restricting the set of particles. For $B_n$, the same
observation applies, and restriction to descendants of $w_m$ preserves
the lexicographic tie convention. Subtracting $V(w_m)$ and deleting the
common genealogical prefix preserve both events. This restriction
property, rather than independence of the events $E_n$, is the basis of
all estimates below.

\subsection{Joint probabilities and tail bounds}

\begin{lemma}[Joint probabilities and exponential moments]
\label{lem:generic-moments}
For $1\leq n_1<\cdots<n_j$, with $n_0=0$,
\begin{equation}\label{eq:generic-multipoint}
 \Q(E_{n_1}\cap\cdots\cap E_{n_j})
 \leq\prod_{i=1}^j r_{n_i-n_{i-1}}.
\end{equation}
Consequently, for every integer $j\geq1$,
\begin{equation}\label{eq:generic-factorial}
 \EQ\left[\binom{U_N}{j}\right]\leq t_N^j.
\end{equation}
For every $\lambda>0$ such that $t_N(e^\lambda-1)<1$,
\begin{equation}\label{eq:generic-mgf}
 \EQ[e^{\lambda U_N}]
 \leq\frac{1}{1-t_N(e^\lambda-1)}.
\end{equation}
\end{lemma}

\begin{proof}
We first remove the latest event from the intersection. The first
$j-1$ events are measurable with respect to $\cH_{n_{j-1}}$, whereas
$E^{(n_{j-1})}_{n_j-n_{j-1}}$ is independent of this sigma-field.
By~\eqref{eq:generic-hereditary}, replacing $E_{n_j}$ by this local
event can only increase the probability. The resulting bound contributes
the factor $r_{n_j-n_{j-1}}$. Iterating proves
\eqref{eq:generic-multipoint}.
Summing over ordered indices and setting $d_i=n_i-n_{i-1}$ gives
\[
 \EQ\left[\binom{U_N}{j}\right]
 \leq\sum_{\substack{d_1,\ldots,d_j\geq1\\d_1+\cdots+d_j\leq N}}
       \prod_{i=1}^jr_{d_i}
 \leq\prod_{i=1}^j\sum_{d_i=1}^Nr_{d_i}=t_N^j.
\]
Finally, expand $(1+z)^{U_N}$ with $z=e^\lambda-1\geq0$.
Taking expectations and bounding the resulting finite sum by the
geometric series proves~\eqref{eq:generic-mgf}.
\end{proof}

\begin{lemma}[A stopping-time lower-tail bound]
\label{lem:generic-lower-tail}
For every $N\geq1$ and every integer $k\geq1$,
\begin{equation}\label{eq:generic-lower-tail}
 \Q(U_N<k)\leq\frac{k}{1+t_N}.
\end{equation}
\end{lemma}

\begin{proof}
Let $T=\inf\{n\geq1:U_n\geq k\}$ be the time of the $k$th
occurrence, with $T=\infty$ if there are fewer than $k$ occurrences.
Since $E_n\in\cH_n$, this is a stopping time for $(\cH_n)$.
After $T$, the events remain dependent, but their expected number is
bounded by that of a fresh descendant process. Indeed, for $n>m$,
\eqref{eq:generic-hereditary} and independence from $\cH_m$ give
\[
 \Q(T=m,E_n)\leq\Q(T=m)r_{n-m}.
\]
Summing over $n=m+1,\ldots,N$ bounds the expected contribution
after $m$, restricted to $\{T=m\}$, by $\Q(T=m)t_{N-m}$.
Before the $k$th occurrence the count is at most $k-1$, and at $T$ it
is exactly $k$. Splitting according to $T$ therefore gives
\begin{align*}
 t_N=\EQ[U_N]
 &\leq(k-1)\Q(T>N)
       +\sum_{m=1}^N\Q(T=m)(k+t_{N-m})\\
 &\leq k-1+(1+t_N)\Q(T\leq N).
\end{align*}
Rearranging and using $\{T>N\}=\{U_N<k\}$ gives the claim.
\end{proof}

\begin{remark}[A direct recurrence criterion]
If $t_N\to\infty$, then for every fixed $k$,
$\Q(\sup_NU_N<k)=0$ by~\eqref{eq:generic-lower-tail}.
Thus infinitely many $E_n$ occur almost surely. If instead
$\sum_nr_n<\infty$, the first Borel--Cantelli lemma gives only
finitely many occurrences almost surely. This argument requires no
zero--one law.
\end{remark}

\subsection{Almost-sure bounds}

The preceding tail estimates do not have to be summable over every
integer $N$. It is enough to apply them when the expected count crosses
a geometric sequence, and then use monotonicity between successive
crossings.

\begin{theorem}[Pathwise bounds in the expected time scale]
\label{thm:generic-pathwise}
Suppose that $t_N\to\infty$. For every $\varepsilon>0$ there is a
deterministic constant $\kappa_\varepsilon>0$ such that almost surely
\begin{equation}\label{eq:generic-pathwise-lower}
 U_N\geq\kappa_\varepsilon
 \frac{t_N}{\log t_N\,(\log\log t_N)^{1+\varepsilon}}
 \qquad\hbox{for all sufficiently large }N.
\end{equation}
Moreover,
\begin{equation}\label{eq:generic-pathwise-upper}
 \limsup_{N\to\infty}\frac{U_N}{t_N\log\log t_N}\leq1
 \qquad\Q\hbox{-almost surely},
\end{equation}
and therefore
\begin{equation}\label{eq:generic-pathwise-exponent}
 \lim_{N\to\infty}\frac{\log U_N}{\log t_N}=1
 \qquad\Q\hbox{-almost surely}.
\end{equation}
All iterated logarithms are used only for sufficiently large $N$.
\end{theorem}

\begin{proof}
\emph{The time grid.}
Fix $\rho>1$ and define deterministic indices
\[
 N_j=\min\{N\geq1:t_N\geq\rho^j\}.
\]
Since $0\leq t_N-t_{N-1}=r_N\leq1$, these indices satisfy
\begin{equation}\label{eq:generic-grid-overshoot}
 \rho^j\leq t_{N_j}\leq\rho^j+1.
\end{equation}
They are finite and strictly increasing for all sufficiently large $j$.
Whenever $N_j\leq N<N_{j+1}$, one has
\begin{equation}\label{eq:generic-grid-interpolation}
 \rho^j\leq t_N\leq\rho^{j+1}+1,\qquad
 \log t_N\sim j\log\rho,\qquad
 \log\log t_N=\log j+O(1),
\end{equation}
uniformly over this interval as $j\to\infty$.

\emph{The lower bound.}
Take $k_j=\lceil t_{N_j}/[j(\log j)^{1+\varepsilon}]\rceil$ for $j\geq2$.
Lemma~\ref{lem:generic-lower-tail} gives
\[
 \Q(U_{N_j}<k_j)
 \leq\frac{1}{j(\log j)^{1+\varepsilon}}+\frac{1}{1+t_{N_j}}.
\]
These probabilities are summable. Therefore almost surely
$U_{N_j}\geq t_{N_j}/[j(\log j)^{1+\varepsilon}]$ for all sufficiently large $j$.
For $N_j\leq N<N_{j+1}$, monotonicity gives
$U_N\geq U_{N_j}$. Equations~\eqref{eq:generic-grid-overshoot}
and~\eqref{eq:generic-grid-interpolation} then prove
\eqref{eq:generic-pathwise-lower}, with a deterministic positive constant
after fixing $\rho=2$.

\emph{The upper bound.}
Fix $0<\theta<1$. Since
$t_N(e^{\theta/t_N}-1)\to\theta$, Lemma~\ref{lem:generic-moments}
and Markov's inequality give constants $K_\theta<\infty$ and
$M_\theta<\infty$ such that
\[
 \Q(U_N\geq xt_N)\leq K_\theta e^{-\theta x}
 \qquad(N\geq M_\theta,\ x\geq0).
\]
For $h>1/\theta$, the probabilities
$\Q(U_{N_j}\geq ht_{N_j}\log j)$ are summable.
Consequently, almost surely $U_{N_j}<ht_{N_j}\log j$ for all
sufficiently large $j$. For $N_j\leq N<N_{j+1}$,
\[
 U_N\leq U_{N_{j+1}}
 \leq h t_{N_{j+1}}\log(j+1).
\]
Using~\eqref{eq:generic-grid-overshoot}
and~\eqref{eq:generic-grid-interpolation}, the limsup in
\eqref{eq:generic-pathwise-upper} is at most $h\rho$.
Let $\theta\uparrow1$, $h\downarrow1$, and $\rho\downarrow1$
along countable sequences satisfying $h>1/\theta$ to obtain the
stated constant $1$.

\emph{The occupation exponent.}
The lower bound with any one fixed $\varepsilon>0$ and the upper
bound imply
$\log U_N=\log t_N+O(\log\log t_N)$ almost surely for all sufficiently
large $N$. Dividing by $\log t_N$ proves
\eqref{eq:generic-pathwise-exponent}.
\end{proof}

\subsection{Tightness at the mean scale}

The almost-sure bounds identify the exponent of $t_N$, but do not
establish a strong law $U_N/t_N\to1$. We now show that $t_N$ is also
the correct scale in probability: the normalized counts are tight and
do not approach zero. A further tail argument shows that any limit in
probability would necessarily be the constant one.

\begin{proposition}[Bounds at the mean scale]
\label{prop:mean-scale-tightness}
Suppose $t_N\to\infty$, and put $X_N=U_N/t_N$.
Then $(X_N)$ is tight and bounded away from zero in probability.
More precisely,
\begin{align}
 \limsup_{N\to\infty}\Q(X_N\leq x)&\leq x
       &&(x\geq0),\label{eq:mean-scale-small}\\
 \limsup_{N\to\infty}\Q(X_N\geq x)&\leq x e^{1-x}
       &&(x>1).\label{eq:mean-scale-large}
\end{align}
Every subsequential weak limit $Y$, with distribution function $F_Y$,
satisfies (with expectations under its limiting law)
\[
 F_Y(x)\leq\min\{x,1\}\quad(x\geq0),\qquad
 \mathbb E[Y]=1,\qquad \mathbb E[Y^j]\leq j!\quad(j\geq1),
\]
and
\[
 \mathbb E[e^{\theta Y}]\leq\frac{1}{1-\theta}\qquad(0\leq\theta<1).
\]
In particular $Y>0$ almost surely. If $X_N$ converges in probability,
its limit is one almost surely.
\end{proposition}

\begin{proof}
\emph{Tail bounds.}
For $x\geq0$, apply~\eqref{eq:generic-lower-tail} with
$k=\lfloor xt_N\rfloor+1$ to obtain
\[
 \Q(X_N\leq x)\leq\frac{xt_N+1}{1+t_N}.
\]
For $0<\theta<1$, equation~\eqref{eq:generic-mgf} gives
\[
 \limsup_N\EQ[e^{\theta X_N}]\leq\frac{1}{1-\theta}.
\]
Markov's inequality therefore bounds the upper-tail limsup by
$e^{-\theta x}/(1-\theta)$. Choosing $\theta=1-1/x$ for $x>1$
proves~\eqref{eq:mean-scale-large}. The two estimates prove the stated
tightness properties.

\emph{Subsequential limits.}
Suppose $X_{N_i}\Rightarrow Y$. At continuity points of the distribution
function of $Y$, the small-tail estimate passes to the limit. Approaching
an arbitrary $x\geq0$ from above through continuity points proves
$F_Y(x)\leq\min\{x,1\}$. For any fixed $0<\theta<1$, the
preceding exponential bound is uniform for all sufficiently large $N$.
Thus $(X_N^j)$ is uniformly integrable for every fixed integer $j$.
In particular $\mathbb E[Y]=\lim_i\EQ[X_{N_i}]=1$.
Expanding ordinary powers in falling factorials and applying
\eqref{eq:generic-factorial} gives
\[
 \limsup_N\EQ[X_N^j]\leq j!,
\]
because the lower-degree terms vanish after division by $t_N^j$.
Uniform integrability passes this bound to $Y$. The exponential bound
for $Y$ follows by applying weak convergence to
$\min\{e^{\theta x},M\}$ and then letting $M\to\infty$.

\emph{Convergence in probability.}
Suppose $X_N\to X$ in probability. Deleting any fixed initial segment
of the spine changes only finitely many event indicators, so it cannot
change this limit. To make this precise, fix $m\geq1$ and let
$U^{(m)}_{N-m}=\sum_{k=1}^{N-m}\1_{E_k^{(m)}}$.
Almost surely all bushes attached before $w_m$ eventually die out.
Beyond that random time, the global and shifted local event indicators
agree, for either choice $E=A$ or $E=B$. Hence
\[
 |U_N-U^{(m)}_{N-m}|=O_m(1)\quad\hbox{almost surely as }N\to\infty,
\]
where the bound is a finite random variable. It follows that
$U^{(m)}_{N-m}/t_N\to X$ in probability. Each variable in this last
sequence is measurable with respect to $\sigma(\mathcal D_i:i\geq m)$.
Taking an almost-surely convergent subsequence shows that $X$ is
measurable with respect to the completion of that sigma-field.
This holds for every $m$, so Kolmogorov's zero--one law makes $X$
almost surely constant. Uniform integrability and $\EQ[X_N]=1$ identify
the constant as one.
\end{proof}

The next corollary combines this general estimate with the survival
asymptotic and, for spatial minima with ties, the multiplicity bound
proved in Section~\ref{sec:multiplicity-bound}.

\begin{corollary}[Logarithmic order in probability]
\label{cor:probability-log-order}
Under finite positive offspring variance, $C_N$ has order $\log N$
in probability for every displacement law. If the displacement law is
nonconstant, the same holds for $D_N$. Explicitly, for every
$\varepsilon>0$ there are deterministic $0<b<B<\infty$ such that
\[
 \liminf_{N\to\infty}
 \Q\bigl(b\log N\leq D_N\leq B\log N\bigr)\geq1-\varepsilon,
\]
and analogously for $C_N$.
\end{corollary}
\begin{proof}
For $C_N$, use $s_N\sim(2/\sigma^2)\log N$. For $D_N$, use
$t_N^A\asymp\log N$, supplied by
Theorem~\ref{thm:uniform-multiplicity}. Apply
Proposition~\ref{prop:mean-scale-tightness} in each case.
\end{proof}

\section{Uniform control of ties at the spatial minimum}
\label{sec:multiplicity-bound}

This section assumes the finite-variance condition \eqref{eq:variance}.
The displacements remain independent across edges and independent of the
family tree. A displacement law is called nonconstant if it is not a
point mass. We show that ties change the order of the one-time
probability by at most a constant factor, even for lattice displacements.

\begin{theorem}[Uniform moments of extremal multiplicity]
\label{thm:uniform-multiplicity}
Let $\mu$ be a nonconstant probability law on $\R$. There is a finite
constant $K=K(p,\mu)$ such that, for every $n\geq1$,
\[
 \E[K_n^2\mid Z_n>0]\leq K.
\]
Consequently, with $a_n=\Q(A_n)$ and
$D_N=\sum_{n=1}^N\1_{A_n}$,
\[
 q_n\leq a_n\leq Kq_n,
 \qquad
 \EQ[D_N]\asymp\log(N+1).
\]
In particular, $a_n\asymp n^{-1}$. The constants may depend on both the
offspring law and the displacement law. This conclusion does not assert
that $n a_n$ has a limit, or that $D_N/\log N$ has an almost-sure limit.
\end{theorem}

The proof has two parts. We first show that the distribution of the
minimum has summably small atoms. Random displacements along the common
ancestral stem provide this anti-concentration even for lattice laws.
This controls the location of the minimum; it does not remove ties
between particles. We then compare the recursion for that location
with recursions weighted by the number of minimizers. The discrepancy
at each generation is controlled by the largest atom, so the total
multiplicative error remains bounded.

For a real random variable $Y$, write
$\operatorname{at}(Y)=\sup_x\Pp(Y=x)$.

\begin{lemma}[Anti-concentration along the common stem]
\label{lem:stem-concentration}
For the minimum $M_n$ defined in the model,
\[
 d_n:=\sup_{x\in\R}\Pp(M_n=x)\leq C(n+1)^{-3/2}.
\]
\end{lemma}

\begin{proof}
\emph{Anti-concentration for a single path.}
If $S_t$ is a sum of $t$ independent variables with the nonconstant law
$\mu$, then
\begin{equation}\label{eq:atom-concentration}
 \operatorname{at}(S_t)\leq C(t+1)^{-1/2}.
\end{equation}
Here is an argument which requires no moment assumption on $\mu$.
If the total mass $r$ of the atomic part of $\mu$ is less than one, any
convolution involving a nonatomic component is nonatomic, so that
$\operatorname{at}(S_t)\leq r^t$. Otherwise choose two distinct atoms
$a,b$ of $\mu$. Condition on which increments belong to $\{a,b\}$ and
on all remaining increments. If $J$ increments belong to this pair,
the conditional sum is a translate and rescaling of a binomial random
variable with $J$ trials and a fixed parameter in $(0,1)$. Its largest
atom is at most $C(J+1)^{-1/2}$. Since $J$ is binomial with $t$ trials
and a fixed positive parameter, averaging gives \eqref{eq:atom-concentration}.

\emph{The length of the common ancestral stem.}
On survival to generation $n$, let $T_n$ be the depth of the most recent
common ancestor of all generation-$n$ vertices; set $T_n=n$ if there is
only one such vertex. For $k\leq n$, let $R_{k,n}$ count vertices at
level $k$ with at least one descendant at level $n$. Then
\[
 \E[R_{k,n}]=q_{n-k},\qquad
 \Pp(R_{k,n}\geq1)=q_n,\qquad
 \{T_n<k,Z_n>0\}=\{R_{k,n}\geq2\}.
\]
It follows that
\[
 \Pp(T_n<k,Z_n>0)\leq q_{n-k}-q_n.
\]
Writing $f(s)=\sum_jp_js^j$, finite variance gives
\[
 q_j-q_{j+1}=f(1-q_j)-(1-q_j)
 \leq \tfrac12\sigma^2q_j^2.
\]
Together with $q_j\asymp(j+1)^{-1}$, this yields
\[
 \Pp(T_n\leq k\mid Z_n>0)
 \leq C\frac{k+1}{n+1},\qquad 0\leq k\leq n.
\]
For $k\leq n/2$ this follows by summing the preceding bound on
$q_j-q_{j+1}$, with $k+1$ in place of $k$; for larger $k$ it follows by
enlarging $C$ and using the trivial bound one.

\emph{Averaging over the stem length.}
Conditional on the finite family tree, every generation-$n$ position
contains the same $T_n$ increments along this common stem. Conditional
also on all increments below the common ancestor, the minimum is a
fixed translate of $S_{T_n}$. Thus
\[
 d_n\leq Cq_n\E[(T_n+1)^{-1/2}\mid Z_n>0].
\]
For $g(k)=(k+1)^{-1/2}$, summation by parts and the bound on the
distribution of $T_n$ give
\begin{align*}
 \E[g(T_n)\mid Z_n>0]
 &=g(n)+\sum_{k=0}^{n-1}(g(k)-g(k+1))
           \Pp(T_n\leq k\mid Z_n>0)\\
 &\leq C(n+1)^{-1/2}.
\end{align*}
Multiplying by $q_n\leq C(n+1)^{-1}$ proves the result.
\end{proof}

\begin{proof}[Proof of Theorem~\ref{thm:uniform-multiplicity}]
\emph{Location and multiplicity measures.}
It is useful to keep track of measures, because the minimum may have a
continuous component even when there is a positive probability of
ties among its particles. Define finite measures on $\R$ by
\[
 \nu_n(E)=\Pp(M_n\in E),\qquad
 \eta_n(E)=\E[K_n\1_{\{M_n\in E\}}],\qquad
 \zeta_n(E)=\E[K_n(K_n-1)\1_{\{M_n\in E\}}].
\]
Each measure places its mass at the minimum location. The weights are,
respectively, one per surviving generation, one per minimizing particle,
and one per ordered pair of distinct minimizing particles. In particular,
\[
 \nu_n(\R)=q_n,\qquad
 \eta_n(\R)=\E[K_n],\qquad
 \zeta_n(\R)=\E[K_n(K_n-1)].
\]
Our aim is to dominate the latter two measures by constant multiples of
$\nu_n$, uniformly in $n$.

Put $\rho_n=\mu*\nu_n$ and
\[
 \overline F_n(x)=1-\rho_n(( -\infty,x]),\qquad
 \overline F_n(x-)=1-\rho_n(( -\infty,x)).
\]
Here $\overline F_n(x-)$ denotes the left limit of the survival
function $\overline F_n$. It is the probability that a child subtree, after its initial
displacement, has no generation-$n$ particle strictly below $x$;
similarly, $\overline F_n(x)$ is the probability that it has no such particle at
or below $x$. Both probabilities include extinction. Their difference
is the mass of the displaced subtree minimum at $x$. Define
\[
 h_n(x)=
 \begin{cases}
 \dfrac{f(\overline F_n(x-))-f(\overline F_n(x))}{\overline F_n(x-)-\overline F_n(x)},&\rho_n(\{x\})>0,\\[6pt]
 f'(\overline F_n(x-)),&\rho_n(\{x\})=0.
 \end{cases}
\]
The minimum-location measure satisfies
\begin{equation}\label{eq:minimum-measure-recursion}
 \nu_{n+1}(dx)=h_n(x)\rho_n(dx).
\end{equation}
At an atom this is the identity
$\nu_{n+1}(\{x\})=f(\overline F_n(x-))-f(\overline F_n(x))$.
For the nonatomic part, let $J$ be the countable atom set of $\rho_n$.
Independent child-subtree minima almost surely cannot tie outside $J$.
For a Borel set $E\subseteq\R\setminus J$, conditioning on $k$ root
children therefore gives
\[
 \Pp(M_{n+1}\in E\mid Z_1=k)
 =k\int_E \overline F_n(x-)^{k-1}\rho_n(dx),\qquad k\geq1.
\]
The contribution for $k=0$ is zero. Summing over $k$ proves
\eqref{eq:minimum-measure-recursion} off $J$ as well.

\emph{Recursions for the multiplicity weights.}
Conditioning on the number of root children and on their independent
displaced subtrees gives two further identities:
\begin{align}
 \eta_{n+1}(dx)
 &=f'(\overline F_n(x-))(\mu*\eta_n)(dx),\label{eq:eta-recursion}\\
 \zeta_{n+1}(dx)
 &=f'(\overline F_n(x-))(\mu*\zeta_n)(dx)\notag\\
 &\quad+f''(\overline F_n(x-))(\mu*\eta_n)(\{x\})
                    (\mu*\eta_n)(dx).\label{eq:zeta-recursion}
\end{align}
For the first identity, a designated child subtree contributes its
minimum multiplicity at $x$, while all other child subtrees have minima
at least $x$. For the second identity, ordered pairs of minimizers either
belong to one child subtree or to two distinct child subtrees. The latter
term is supported on atoms of $\mu*\eta_n$, as displayed.

\emph{A summable comparison error.}
The minimum-location recursion involves the secant slope $h_n$, whereas
the multiplicity recursions involve the endpoint derivative $f'(\overline F_n(x-))$.
These agree away from atoms. At atoms their difference is small enough
to sum over generations.
Since the critical offspring law is nondegenerate, $p_0>0$ and
$c_0:=f'(p_0)>0$. For $n\geq1$, both $\overline F_n(x-)$ and $\overline F_n(x)$ are at
least $\Pp(Z_n=0)\geq p_0$. Also $f''(s)\leq\sigma^2$ on $[0,1]$.
Consequently,
\begin{equation}\label{eq:secant-control}
 h_n(x)\geq c_0,\qquad
 f'(\overline F_n(x-))\leq h_n(x)(1+C d_n).
\end{equation}
Indeed, the discrepancy between the endpoint derivative and its secant
average is at most
$\sigma^2\rho_n(\{x\})\leq\sigma^2d_n$; for a nonatom it is zero.

\emph{The first multiplicity moment.}
The lower bound $\overline F_n(x-),\overline F_n(x)\geq p_0$ applies from generation one onward,
so we treat the initial comparison separately.
For the first generation one still has $\eta_1\leq \kappa_1\nu_1$ with
$\kappa_1<\infty$. To see this, convexity gives, for $0\leq v<u\leq1$,
\[
 \frac{f'(u)}{(f(u)-f(v))/(u-v)}
 \leq \frac{u f'(u)}{f(u)-p_0}.
\]
The right-hand side is bounded on $(0,1]$: it tends at zero to the
smallest positive $j$ with $p_j>0$, and it is continuous elsewhere.
At nonatomic points the ratio is one.
Since $\eta_0=\nu_0=\delta_0$, the assertion follows from
\eqref{eq:minimum-measure-recursion}--\eqref{eq:eta-recursion}.

Now \eqref{eq:secant-control} shows inductively that
\[
 \eta_n\leq \kappa_1\prod_{j=1}^{n-1}(1+C d_j)\,\nu_n.
\]
Lemma~\ref{lem:stem-concentration} implies $\sum_jd_j<\infty$, so
$\eta_n\leq \kappa_2\nu_n$ uniformly in $n$.

\emph{The second multiplicity moment.}
The first-generation second factorial measure also
satisfies $\zeta_1\leq \kappa_3\nu_1$. At an atom $x$ of $\mu$, let
$u=\mu([x,\infty))$ and $v=\mu((x,\infty))$. Its ratio is
\[
 \frac{f''(u)(u-v)^2}{f(u)-f(v)}
 \leq \frac{u^2 f''(u)}{f(u)-p_0}.
\]
This is bounded on $(0,1]$ by the same small-$u$ argument and finite
offspring variance. The nonatomic part of $\zeta_1$ is zero.
Using $\eta_n\leq \kappa_2\nu_n$, the forcing term in
\eqref{eq:zeta-recursion} is bounded by
$\sigma^2\kappa_2^2d_n\rho_n(dx)\leq C d_n\nu_{n+1}(dx)$.
Thus, if $\zeta_n\leq b_n\nu_n$, then
\[
 \zeta_{n+1}\leq\bigl((1+C d_n)b_n+C'd_n\bigr)\nu_{n+1}.
\]
Summability of $d_n$ proves that $b_n$ remains bounded: both the
multiplicative factors and the accumulated additive terms in this
recursion are uniformly bounded. Taking total
masses yields $\E[K_n^2]\leq Kq_n$. Finally
$K_n\geq\1_{\{Z_n>0\}}$ and
$\Q(A_n)=\E[K_n]$, proving all claims.
\end{proof}

\begin{remark}[The constant-displacement exception]
If $\mu=\delta_c$, every generation-$n$ particle lies at $nc$, so
$A_n$ holds at every generation under $\Q$. Thus nonconstancy is
necessary for the preceding $n^{-1}$ order and logarithmic expected
frequency. For atomless $\mu$ the sharper identity $a_n=q_n$ holds.
\end{remark}

\section{Recurrence as the unique spatial minimum}
\label{sec:strict-minimum}

Define the strict-minimum event by
\[
 A_n^{\mathrm{str}}=
 \{V(w_n)<V(u)\text{ for every }u\in\cT_n\setminus\{w_n\}\},
 \qquad D_N^{\mathrm{str}}=\sum_{n=1}^N\1_{A_n^{\mathrm{str}}}.
\]
The inequality is vacuous if the generation has only one particle.

\begin{theorem}[Strict spatial extremality]
\label{thm:strict-minimum}
Assume \eqref{eq:critical} and \eqref{eq:variance}.
If $\mu$ is nonconstant, there exists $c=c(p,\mu)>0$ such that
\[
 c q_n\leq\Q(A_n^{\mathrm{str}})\leq q_n
 \qquad(n\geq1).
\]
Consequently, the spine is the unique spatial minimum infinitely often
almost surely, and
\[
 \EQ[D_N^{\mathrm{str}}]\asymp\log(N+1),\qquad
 \frac{\log D_N^{\mathrm{str}}}{\log\log N}\longrightarrow1
 \quad\Q\text{-almost surely}.
\]
Moreover, $D_N^{\mathrm{str}}$ has order $\log N$ in probability.
If $\mu$ is a point mass, only finitely many strict-minimum events occur
almost surely. In fact
\[
 \EQ\left[\sum_{n=0}^{\infty}\1_{A_n^{\mathrm{str}}}\right]\leq\frac1{p_0}.
\]
\end{theorem}

\begin{proof}
\emph{The recursion for a unique minimizer.}
The finite-generation spine identity gives
\begin{equation}\label{eq:strict-identity}
 \Q(A_n^{\mathrm{str}})=\Pp(K_n=1).
\end{equation}
Define the subprobability measure
$\xi_n(dx)=\Pp(K_n=1,M_n\in dx)$.
We retain $\nu_n,\rho_n,\overline F_n,h_n$ from the proof of
Theorem~\ref{thm:uniform-multiplicity}. A unique overall minimizer
must be the unique minimizer of one root-child subtree, while the minima
of all other child subtrees must be strictly larger. Thus
\begin{equation}\label{eq:strict-measure-recursion}
 \xi_{n+1}(dx)=f'(\overline F_n(x))(\mu*\xi_n)(dx).
\end{equation}
For a finite positive measure $\alpha$, write
\[
 (\mathcal K_n\alpha)(dx)=f'(\overline F_n(x))(\mu*\alpha)(dx).
\]
Thus each operator first adds an independent displacement, then weights
the resulting position by $f'(\overline F_n(x))$. This factor sums over
the choice of a distinguished child and requires every competing child
subtree to lie strictly to its right. The unique-minimum measure
follows these operators exactly. The ordinary minimum measure follows
the same recursion up to a summable multiplicative error.
For $n\geq1$, the secant comparison used earlier also gives
\[
 f'(\overline F_n(x))\leq h_n(x)
 \leq f'(\overline F_n(x))(1+C d_n).
\]
Therefore, with
$P_*=\prod_{j=1}^{\infty}(1+C d_j)<\infty$,
\begin{equation}\label{eq:strict-semigroup-comparison}
 \nu_n\leq P_* \mathcal K_{n-1}\cdots \mathcal K_1\nu_1,
 \qquad
 \xi_n=\mathcal K_{n-1}\cdots \mathcal K_1\xi_1.
\end{equation}
For $n=1$, the empty composition is the identity.

\emph{Initial mass on a lower halfline.}
It remains to compare the initial measures on positions from which the
operators retain enough mass. Because $\mu$ is nonconstant, choose
$b\in\R$ such that
$\mu(( -\infty,b])>0$ and $r:=\mu((b,\infty))>0$.
Put $E=( -\infty,b]$. For $x\in E$, one has
$\overline F_0(x)=\mu((x,\infty))\geq r$. Since $h_0(x)\leq f'(1)=1$,
\[
 \xi_1|_E\geq f'(r)\mu|_E\geq f'(r)\nu_1|_E.
\]
Moreover $\nu_1(E)>0$ and $f'(r)>0$.

\emph{Propagation of the comparison.}
Let $g_n(x)$ be the total mass of
$\mathcal K_{n-1}\cdots \mathcal K_1\delta_x$. Each function
$x\mapsto f'(\overline F_j(x))$ is nonincreasing. Coupling the displacements
in the two operator products started from $x\leq y$ consequently gives
$g_n(x)\geq g_n(y)$. This remains true without any moment assumption
on the displacements: with identical increments, the path started
farther left receives at least as large a weight at every step.
For a nonincreasing nonnegative function $g_n$,
\[
 \int_E g_n\,d\nu_1
 \geq\frac{\nu_1(E)}{\nu_1(\R)}\int_\R g_n\,d\nu_1.
\]
Indeed, the average of $g_n$ over a lower halfline is at least its
average over the complementary upper halfline. Taking total masses
in \eqref{eq:strict-semigroup-comparison} now yields
\[
 \Pp(K_n=1)=\xi_n(\R)
 \geq \frac{f'(r)\nu_1(E)}{q_1 P_*}\,q_n.
\]
The upper bound $\Pp(K_n=1)\leq q_n$ is immediate.

\emph{Recurrence and frequency.}
The strict-minimum events have the same hereditary inclusion and
independent descendant-tree structure as the weak-minimum events:
a unique minimum remains unique after restriction to a subset containing
it. Hence the arguments of Section~\ref{sec:generic-frequency} apply
verbatim with $E_n=A_n^{\mathrm{str}}$. The probability comparison
and $q_n\asymp n^{-1}$ give the recurrence and occupation claims.

\emph{The constant-displacement case.}
Finally, if $\mu=\delta_a$, every particle at generation $n$ has
position $na$, so $A_n^{\mathrm{str}}=\{Z_n=1\}$ under $\Q$.
Size bias gives $\Q(Z_n=1)=\Pp(Z_n=1)$.
Under $\Pp$, on each visit of the population process to state one,
the next step is extinction with probability $p_0>0$.
By the strong Markov property, the total number of visits to state one
(including the initial visit) is stochastically dominated by a geometric
random variable with success probability $p_0$. Thus
\[
 \sum_{n\geq0}\Q(Z_n=1)
 =\sum_{n\geq0}\Pp(Z_n=1)\leq p_0^{-1}.
\]
The first Borel--Cantelli lemma proves the final assertion.
\end{proof}

\section{Consequences for spatial extremal generations}
\label{sec:consequences}

\begin{proof}[Proof of Theorem~\ref{thm:frequency-main}]
For the selected event, $r_n=q_n$ and $t_N=s_N\sim c\log N$, where
$c=2/\sigma^2$. Theorem~\ref{thm:generic-pathwise} therefore gives
\eqref{eq:selected-exponent}.
For nonconstant $\mu$, Theorem~\ref{thm:uniform-multiplicity} gives
$q_n\leq a_n\leq Kq_n$, so $t_N^A\asymp\log N$ and
$\log t_N^A/\log\log N\to1$. Applying
Theorem~\ref{thm:generic-pathwise} to $E_n=A_n$ proves
\eqref{eq:weak-comparison}. In particular, $D_N/N\to0$. The assertions about probability order
follow from Corollary~\ref{cor:probability-log-order}.
The direct recurrence criterion gives $\tau_k<\infty$ almost surely.
Since the count increases by one at every occurrence,
$D_{\tau_k}=k$ and $\tau_k\to\infty$. Substitution in
\eqref{eq:weak-comparison} gives
$\log k/\log\log\tau_k\to1$, which is equivalent to
\eqref{eq:arrival-exponent}. The constant-displacement assertion follows
because all generation-$n$ positions then equal $nc$.
\end{proof}

\begin{corollary}[Explicit bounds for the selected count]
\label{cor:explicit-count}
Under finite positive offspring variance, put $c=2/\sigma^2$.
For every $\varepsilon>0$, there is a deterministic $b_\varepsilon>0$
such that almost surely, for all sufficiently large $N$,
\[
 C_N\geq b_\varepsilon
 \frac{\log N}{\log\log N\,(\log\log\log N)^{1+\varepsilon}}.
\]
Also
\[
 \limsup_{N\to\infty}
 \frac{C_N}{c\log N\log\log\log N}\leq1
 \qquad\Q\text{-almost surely}.
\]
For nonconstant $\mu$, the same orders of lower and upper bounds hold
for $D_N$, with constants depending on $p$ and $\mu$.
\end{corollary}
\begin{proof}
Substitute $s_N\sim c\log N$ in the bounds of
Theorem~\ref{thm:generic-pathwise}. For $D_N$ use
$t_N^A\asymp\log N$ instead. All expressions are considered only for
sufficiently large $N$.
\end{proof}

\begin{proposition}[A two-time lower bound]
For every $m,k\geq1$,
\[
 \frac{q_k}{1-q_k}(q_m-q_{m+k})
 \leq \Q(B_m\cap B_{m+k})\leq q_mq_k.
\]
\end{proposition}
\begin{proof}
Let $\mathcal G_m$ record the marked spatial tree through generation $m$,
but not reproduction at generation $m$. Given $\mathcal G_m$, the
descendant tree of $w_m$ has law $\Q$, and those of the other $Z_m-1$
vertices have law $\Pp$; these descendant trees are independent. If all
the latter trees are extinct at relative generation $k$, the whole
generation $m+k$ consists of descendants of $w_m$. Consequently
\[
 \Q(B_m\cap B_{m+k})
 \geq q_k\EQ[\1_{B_m}(1-q_k)^{Z_m-1}].
\]
The finite-generation spine identity, applied to the selected vertex,
turns the right-hand expectation into
$\E[\1_{\{Z_m>0\}}(1-q_k)^{Z_m-1}]$.
Writing $f_m(s)=\E[s^{Z_m}]$ and using the branching identity
$f_m(1-q_k)=1-q_{m+k}$ gives the asserted lower bound. The denominator
is positive: criticality and nondegeneracy imply $p_0>0$, hence $q_k<1$.
The upper bound was proved earlier.
\end{proof}

\begin{proposition}[Ancestral persistence and separated times]
\label{prop:ancestral-persistence}
Write $u|_m$ for the generation-$m$ ancestor of a vertex $u$. For every
$m,k\geq1$,
\[
 \Q(B_m\cap B_{m+k})
 =\Pp\bigl(Z_{m+k}>0,\ \ell_{m+k}|_m=\ell_m\bigr).
\]
Under finite positive offspring variance, along any sequence with
$k/m\to\infty$,
\[
 \Q(B_m\cap B_{m+k})\sim q_mq_k\sim q_mq_{m+k}.
\]
\end{proposition}
\begin{proof}
Apply the finite-generation spine identity at time $m+k$. In the
resulting sum, the event $B_{m+k}$ selects exactly the vertex
$\ell_{m+k}$, while $B_m$ requires its generation-$m$ ancestor to be
$\ell_m$. This proves the identity. For the asymptotics, the two-time
bounds give
\[
 \frac{1-q_{m+k}/q_m}{1-q_k}
 \leq\frac{\Q(B_m\cap B_{m+k})}{q_mq_k}\leq1.
\]
The estimates $q_n\asymp(n+1)^{-1}$ imply $q_{m+k}/q_m\to0$ and
$q_k\to0$ when $k/m\to\infty$. Finally
$q_k/q_{m+k}\to1$ by the survival asymptotic.
\end{proof}

\begin{corollary}[Arrival times in probability]
\label{cor:arrival-probability-scale}
Under finite positive offspring variance, let
$\tau_j^B=\inf\{N:C_N\geq j\}$ and put $c=2/\sigma^2$.
Then $\log\tau_j^B$ has order $j$ in probability. More precisely,
\[
 \limsup_{j\to\infty}\Q\left(\frac{c\log\tau_j^B}{j}>R\right)
 \leq\frac1R\qquad(R>0),
\]
and
\[
 \limsup_{j\to\infty}\Q\left(\frac{c\log\tau_j^B}{j}\leq x\right)
 \leq x^{-1}e^{1-1/x}\qquad(0<x<1).
\]
For nonconstant displacements, the weak-minimum arrival times also
have $\log\tau_j$ of order $j$ in probability.
\end{corollary}
\begin{proof}
For $N_j=\lfloor\exp(Rj/c)\rfloor$,
\[
 \Q(\tau_j^B>N_j)=\Q(C_{N_j}<j)\leq\frac{j}{1+s_{N_j}},
\]
whose upper limit is $1/R$. For
$N_j^-=\lfloor\exp(xj/c)\rfloor$, the exponential-moment bound, with
parameter $\theta/s_{N_j^-}$, gives
\[
 \limsup_j\Q(\tau_j^B\leq N_j^-)
 =\limsup_j\Q(C_{N_j^-}\geq j)
 \leq\frac{e^{-\theta/x}}{1-\theta}.
\]
Take $\theta=1-x$. These bounds prove the probability order for the
selected arrival times. The same argument for $D_N$, using
$t_N^A\asymp\log N$, proves the weak-minimum version with unspecified
positive constants.
\end{proof}

\section{Selected-count fluctuations}\label{sec:fluctuations}

\begin{proposition}[Some pathwise information at the mean scale]
If $s_N\to\infty$, then
\[
 \liminf_N C_N/s_N\leq1\leq\limsup_N C_N/s_N
 \qquad\Q\text{-almost surely}.
\]
Both the liminf and limsup are almost surely deterministic extended
nonnegative constants. If $C_N/s_N$ has an almost-sure limit, that limit
must equal one.
\end{proposition}
\begin{proof}
Deleting the first $m$ spine blocks changes the counting sequence,
after a shift by $m$, by only a finite random amount: every off-spine
tree in those blocks eventually becomes extinct. Moreover,
$s_{N-m}/s_N\to1$, because $0\leq s_N-s_{N-m}\leq m$.
Thus both the liminf and limsup are tail-measurable and are almost surely
constant. Fatou's lemma and $\EQ[C_N]/s_N=1$ bound the liminf by one.
For $0<\eta<1$, the second-moment bound and Paley--Zygmund give
\[
 \liminf_N\Q(C_N/s_N\geq\eta)\geq(1-\eta)^2/2>0.
\]
Therefore the limsup is at least $\eta$ with positive probability,
and hence almost surely by its deterministic character. Let
$\eta\uparrow1$ through a countable sequence.
\end{proof}

\subsection{An exactly solvable tied model}

\begin{theorem}[Deterministic displacements]
Assume $\mu=\delta_a$ and $0<\sigma^2<\infty$. Then
\[
 \frac{C_N}{s_N}\ \xrightarrow{\mathrm d}\ \operatorname{Exp}(1),
 \qquad s_N\sim\frac{2}{\sigma^2}\log N.
\]
Furthermore,
\[
 \liminf_N C_N/s_N=0,
 \qquad\limsup_N C_N/s_N=\infty
 \qquad\Q\text{-almost surely}.
\]
This theorem concerns the lexicographically selected minimum. The count
of weak spatial-minimum events is exactly $N$ in this model.
\end{theorem}
\begin{proof}
All generation-$n$ positions equal $na$, so $B_n$ means that the spine
is the genealogically leftmost particle. On $B_m$, every other vertex
in generation $m$ is lexicographically to the right of $w_m$; all of
its descendants remain to the right of every descendant of $w_m$.
Consequently
\[
 B_m\cap B_{m+k}=B_m\cap B_k^{(m)}.
\]
Independence of the future spine blocks now gives the exact multipoint
identity
\[
 \Q\Bigl(\bigcap_{i=1}^r B_{n_i}\Bigr)
 =\prod_{i=1}^r q_{n_i-n_{i-1}},\qquad 0=n_0<n_1<\cdots<n_r.
\]
It follows that
\[
 s_{\lfloor N/r\rfloor}^{\,r}
 \leq\EQ\left[\binom{C_N}{r}\right]
 =\sum_{\substack{k_1,\ldots,k_r\geq1\\\sum k_i\leq N}}
     \prod_iq_{k_i}
 \leq s_N^r.
\]
For each fixed $r$, $s_{\lfloor N/r\rfloor}/s_N\to1$. Thus every
factorial moment satisfies $\EQ[(C_N)_r]/s_N^r\to r!$, and the lower
factorial moments contribute negligibly when ordinary moments are
expanded. Therefore $\EQ[(C_N/s_N)^r]\to r!$ for every $r$.
The exponential law is moment-determinate, proving convergence in law.

For every $\varepsilon>0$, convergence in law gives
$\Q(C_N/s_N<\varepsilon)\to1-e^{-\varepsilon}>0$, so the liminf is
at most $\varepsilon$ with positive probability and hence almost surely
by tail triviality. Similarly, for each $M>0$,
$\Q(C_N/s_N>M)\to e^{-M}>0$, so the limsup is at least $M$ almost surely.
Take countable sequences $\varepsilon\downarrow0$ and $M\uparrow\infty$.
\end{proof}

\begin{remark}[Limits of these conclusions]
The deterministic-displacement example proves that an almost-sure
logarithmic asymptotic for the selected count is false without further
assumptions. It does not disprove such an asymptotic for a nondegenerate
atomless displacement law. Neither the two-time bounds nor the
multipoint upper bound above establish an exponential limit, or any
particular nondegenerate limit, in that case.
\end{remark}

\section{A finite-occurrence regime with infinite offspring variance}
\label{sec:stable-multiplicity}

\begin{proposition}[Stable offspring and nonconstant displacements]
\label{prop:stable-multiplicity}
Suppose that the offspring generating function is
\[
 f(s)=s+c(1-s)^\alpha,\qquad
 1<\alpha<2,\qquad 0<c\leq1/\alpha,
\]
and that the independent edge-displacement law $\mu$ is not a point mass.
Put $\gamma=\alpha-1$ and $\beta=1/\gamma$.
There is a finite constant $K=K(\alpha,c,\mu)$ such that
\[
 q_n\sim[c(\alpha-1)n]^{-\beta},\qquad
 1\leq\E[K_n\mid Z_n>0]\leq K\quad(n\geq1).
\]
Consequently,
\[
 q_n\leq\Q(A_n)\leq Kq_n,\qquad
 \EQ\left[\sum_{n\geq1}\1_{A_n}\right]<\infty.
\]
In particular, the spine is at the weak spatial minimum in only finitely
many generations almost surely.
\end{proposition}

The proof follows the minimum-location comparison used in
Section~\ref{sec:multiplicity-bound}. The new point is that $f'$ is only
H\"older continuous. The resulting comparison errors are still summable
because the survival probabilities now decay faster than $n^{-1}$.

\begin{proof}
\emph{Survival asymptotics.}
The displayed function is a critical, nondegenerate offspring
generating function: $p_0=c$, $p_1=1-c\alpha\geq0$, and its remaining
coefficients are positive. Its mean is one, whereas its variance is
infinite. The survival recursion is exact:
\[
 q_{n+1}=q_n-cq_n^\alpha.
\]
It gives $q_n\downarrow0$ and
$q_{n+1}^{-\gamma}-q_n^{-\gamma}\to c\gamma$.
Taking arithmetic means proves the asserted asymptotic. In particular,
$q_n\asymp(n+1)^{-\beta}$.

\emph{Atoms of the minimum.}
We next bound the atoms of the minimum location $M_n$.
On survival to generation $n$, let $T_n$ be the depth of the most recent
common ancestor of all generation-$n$ particles, with $T_n=n$ for a
singleton generation. The reduced-tree calculation used in
Lemma~\ref{lem:stem-concentration} does not require finite variance:
\[
 \Pp(T_n<k,Z_n>0)\leq q_{n-k}-q_n.
\]
For $k\leq n/2$, sum the exact increments
$q_j-q_{j+1}=cq_j^\alpha$ and use $\alpha\beta=\beta+1$.
After division by $q_n$, and with $k+1$ in place of $k$, this gives
\[
 \Pp(T_n\leq k\mid Z_n>0)
 \leq C\frac{k+1}{n+1},\qquad 0\leq k\leq n.
\]
The remaining range follows from the trivial bound one, by increasing
$C$. The atom estimate~\eqref{eq:atom-concentration} requires only
nonconstancy of $\mu$. Conditioning on the family tree and on all
increments below the common ancestor therefore gives
\begin{align*}
 d_n:=\sup_{x\in\R}\Pp(M_n=x)
 &\leq Cq_n\E[(T_n+1)^{-1/2}\mid Z_n>0]\\
 &\leq Cq_n(n+1)^{-1/2}.
\end{align*}
The last step is summation by parts using the preceding distributional
bound. Consequently,
\begin{equation}\label{eq:stable-summable-atoms}
 \sum_{n\geq1}d_n^\gamma<\infty,
 \qquad
 d_n^\gamma\leq C(n+1)^{-1-\gamma/2}.
\end{equation}

\emph{Comparison of recursive measures.}
Use the minimum-location and first-multiplicity measures
\[
 \nu_n(E)=\Pp(M_n\in E),\qquad
 \eta_n(E)=\E[K_n\1_{\{M_n\in E\}}].
\]
Their recursions~\eqref{eq:minimum-measure-recursion}
and~\eqref{eq:eta-recursion} require only finite offspring mean.
In the notation of those recursions,
$\rho_n=\mu*\nu_n$,
$\overline F_n(x)=1-\rho_n(( -\infty,x])$, and
$\overline F_n(x-)=1-\rho_n(( -\infty,x))$.
The function $h_n(x)$ is the secant slope of
$f$ on $[\overline F_n(x),\overline F_n(x-)]$, interpreted as $f'(\overline F_n(x-))$ when the
endpoints coincide.

Here $f'(s)=1-c\alpha(1-s)^\gamma$ is H\"older continuous of
order $\gamma$ on $[0,1]$:
\[
 |f'(u)-f'(v)|\leq c\alpha|u-v|^\gamma.
\]
For $n\geq1$, both endpoints are at least
$1-q_n\geq p_0=c$, so $h_n(x)\geq f'(c)>0$. Since
$\rho_n(\{x\})\leq d_n$, the endpoint derivative and secant slope
satisfy
\[
 0\leq f'(\overline F_n(x-))-h_n(x)\leq c\alpha d_n^\gamma,
 \qquad
 f'(\overline F_n(x-))\leq h_n(x)(1+C d_n^\gamma).
\]

\emph{Iteration and finite occurrence.}
The first-generation domination $\eta_1\leq \kappa_1\nu_1$ also needs
only finite mean. Indeed, the derivative-to-secant ratio at an atom
with endpoints $0\leq v<u\leq1$ is bounded by
\[
 \frac{f'(u)}{[f(u)-f(v)]/(u-v)}
 \leq\frac{u f'(u)}{f(u)-p_0}.
\]
The right-hand side is bounded: its limit at zero is the smallest
positive offspring size with positive probability, and its value at
one is finite. At nonatomic points the ratio is one.
Induction in the two measure recursions now yields
\[
 \eta_n\leq \kappa_1\prod_{j=1}^{n-1}(1+C d_j^\gamma)\nu_n
 \leq K\nu_n,
\]
where the last inequality follows from
\eqref{eq:stable-summable-atoms}. Taking total masses gives
$\E[K_n]\leq Kq_n$, while $K_n\geq\1_{\{Z_n>0\}}$ gives the
opposite inequality. Finally use~\eqref{eq:multiplicity} and
$\beta>1$ to obtain the finite expected total count and the
almost-sure conclusion.
\end{proof}

\begin{remark}
This proposition bounds the conditional first moment of multiplicity,
not its second moment. A second-moment analogue is false for this
offspring family: if $\mu$ has a smallest support point with positive
mass $r$, the number of first-generation children at that point has
second factorial moment $r^2\E[Z_1(Z_1-1)]=\infty$ and is bounded above
by $K_1$. Nonconstancy remains necessary for the finite-occurrence
conclusion, since point-mass displacements give $A_n$ in every generation.
\end{remark}

\subsection{Last extremal generations and quantitative separation}
\label{sec:finite-branch-extensions}

\begin{proposition}[The last occurrence has the survival-tail order]
\label{prop:last-occurrence-tail}
Fix either $E_n=A_n$ for every $n$, or $E_n=B_n$ for every $n$.
Write $r_n=\Q(E_n)$, and suppose
$t_\infty=\sum_{n\geq1}r_n<\infty$. Define
\[
 L_E=\sup\{n\geq1:E_n\text{ occurs}\},
 \qquad \sup\varnothing=0,
 \qquad b_N=\sum_{n>N}r_n.
\]
Then, for every integer $N\geq0$,
\begin{equation}\label{eq:last-occurrence-tail}
 \frac{b_N}{1+t_\infty}
 \leq\Q(L_E>N)\leq b_N.
\end{equation}
In particular, if $r_n\asymp n^{-\beta}$ for some $\beta>1$, then
\begin{equation}\label{eq:last-occurrence-power}
 \Q(L_E>N)\asymp N^{1-\beta},\qquad
 \EQ[L_E^u]<\infty\ \Longleftrightarrow\ 0<u<\beta-1
 \quad(u>0).
\end{equation}
\end{proposition}

\begin{proof}
The upper bound is the union bound. For the lower bound, let
$T=\inf\{n>N:E_n\text{ occurs}\}$, with $T=\infty$ if the set is empty.
For $m>N$, the event $\{T=m\}$ belongs to $\cH_m$. The hereditary
inclusion and independence of the future spine blocks give
\[
 \Q(T=m,E_n)\leq\Q(T=m)r_{n-m}\qquad(n>m).
\]
For $M>N$, decompose the number of occurrences in $(N,M]$ according
to its first occurrence. This yields
\begin{align*}
 \sum_{n=N+1}^{M}r_n
 &\leq\sum_{m=N+1}^{M}\Q(T=m)
             \left(1+\sum_{j=1}^{M-m}r_j\right)\\
 &\leq(1+t_\infty)\Q(T\leq M).
\end{align*}
Let $M\to\infty$ to prove the lower bound in
\eqref{eq:last-occurrence-tail}.
If $r_n\asymp n^{-\beta}$, then $b_N\asymp N^{1-\beta}$.
For $u>0$, the identity
\[
 \EQ[L_E^u]
 =\sum_{n\geq0}\bigl((n+1)^u-n^u\bigr)\Q(L_E>n)
\]
shows that the moment is finite exactly when
$\sum_{n\geq1}n^{u-\beta}<\infty$, namely $u<\beta-1$.
\end{proof}

We next give a spatial estimate which complements finiteness of the
extremal-generation set. Put
\[
 G_n=V(w_n)-M_n\geq0.
\]
The key observation is that moving the final edge of the spine particle
changes the location of that particle alone within generation $n$. A shift
of this edge turns a near-minimum event into an exact-minimum event, at the
cost of a likelihood ratio. Combining this identity with the decay of
$q_n$ gives summable bounds for near-minimum events.

\begin{lemma}[A terminal-edge change of measure]
\label{lem:terminal-edge-shift}
Suppose that the displacement law has an everywhere positive density
$g$ on $\R$. For $h\geq0$, put
$R_h(y)=g(y+h)/g(y)$, and let $Y_n=V(w_n)-V(w_{n-1})$.
Then, for every $n\geq1$,
\begin{equation}\label{eq:terminal-edge-identity}
 \Q(G_n\leq h)=\EQ\bigl[R_h(Y_n)\1_{A_n}\bigr].
\end{equation}
Consequently, for every $p\in(1,\infty)$ with
$\int R_h(y)^p g(y)\,dy<\infty$,
\begin{equation}\label{eq:terminal-edge-holder}
 \Q(G_n\leq h)
 \leq\left(\int R_h(y)^p g(y)\,dy\right)^{1/p}
          q_n^{1-1/p}.
\end{equation}
If $R_h$ is essentially bounded with respect to $g(y)\,dy$, then
\begin{equation}\label{eq:terminal-edge-bounded}
 \Q(G_n\leq h)\leq\|R_h\|_\infty q_n.
\end{equation}
\end{lemma}

\begin{proof}
Condition on the genealogical tree, its distinguished spine, and all
edge displacements except the edge from $w_{n-1}$ to $w_n$.
Conditional on this information, $Y_n$ has density $g$.
Let
\[
 W_n=\min\{V(u):u\in\cT_n,\ u\ne w_n\},
 \qquad H_n=W_n-V(w_{n-1}),
\]
where the minimum over the empty set is $+\infty$.
Every other generation-$n$ position is unchanged when $Y_n$ is varied,
so $H_n$ is measurable with respect to the conditioned information.
For finite $H_n$, the events in question are
$A_n=\{Y_n\leq H_n\}$ and
$\{G_n\leq h\}=\{Y_n\leq H_n+h\}$; when $H_n=+\infty$ both events
hold and the same integral identity applies. Changing variables gives
\[
 \int_{-\infty}^{H_n+h}g(y)\,dy
 =\int_{-\infty}^{H_n}g(y+h)\,dy
 =\int_{-\infty}^{H_n}R_h(y)g(y)\,dy.
\]
Taking expectations proves~\eqref{eq:terminal-edge-identity}.
The density assumption implies that $\mu$ is atomless, so
$\Q(A_n)=q_n$. H\"older's inequality gives
\eqref{eq:terminal-edge-holder}; the essential-supremum estimate gives
\eqref{eq:terminal-edge-bounded}.
\end{proof}

\begin{theorem}[Quantitative eventual separation from the minimum]
\label{thm:quantitative-gap}
Suppose that $q_n\leq C n^{-\beta}$ for some $\beta>1$ and finite $C$.
The following assertions hold under the marked size-biased law.
\begin{enumerate}
\item If $\mu$ is Gaussian with arbitrary mean and variance $v>0$, then
\begin{equation}\label{eq:gaussian-gap-lower}
 \liminf_{n\to\infty}\frac{G_n}{\sqrt{\log n}}
 \geq\sqrt{2v}\,(\sqrt\beta-1)
 \qquad\Q\text{-almost surely}.
\end{equation}
\item If $\mu$ is a Laplace law with arbitrary location and scale $b>0$,
having density $g(y)=(2b)^{-1}\exp(-|y-m|/b)$, then
\begin{equation}\label{eq:laplace-gap-lower}
 \liminf_{n\to\infty}\frac{G_n}{\log n}
 \geq b(\beta-1)
 \qquad\Q\text{-almost surely}.
\end{equation}
\end{enumerate}
In particular, $G_n\to\infty$ almost surely in both cases. These are
eventual lower bounds, not assertions of an exact spatial growth rate.
\end{theorem}

\begin{proof}
\emph{Gaussian increments.}
For a Gaussian density with mean $m$ and variance $v$,
\[
 R_h(y)=\exp\left(-\frac{h(y-m)}v-\frac{h^2}{2v}\right),
 \qquad
 \left(\int R_h(y)^pg(y)\,dy\right)^{1/p}
 =\exp\left(\frac{(p-1)h^2}{2v}\right).
\]
Lemma~\ref{lem:terminal-edge-shift} therefore gives
\[
 \Q(G_n\leq h)
 \leq C^{1-1/p}
       \exp\left(\frac{(p-1)h^2}{2v}\right)
       n^{-\beta(1-1/p)}.
\]
Fix $a$ with $0<a<\sqrt{2v}(\sqrt\beta-1)$, and take
$h=a\sqrt{\log n}$ for $n\geq2$. Writing $d=a^2/(2v)$,
choose $p=\sqrt{\beta/d}>1$. The exponent in the resulting bound is
\[
 \beta(1-1/p)-d(p-1)
 = (\sqrt\beta-\sqrt d)^2>1.
\]
Thus $\sum_n\Q(G_n\leq a\sqrt{\log n})<\infty$.
The first Borel--Cantelli lemma and a countable sequence of $a$ tending
up to $\sqrt{2v}(\sqrt\beta-1)$ prove
\eqref{eq:gaussian-gap-lower}.

\emph{Laplace increments.}
For the Laplace density, the triangle inequality gives
\[
 R_h(y)=\exp\left(\frac{|y-m|-|y+h-m|}{b}\right)
 \leq e^{h/b}.
\]
Consequently, for $0<a<b(\beta-1)$,
\[
 \Q(G_n\leq a\log n)\leq C n^{-\beta+a/b},
\]
and the probabilities are summable. Apply Borel--Cantelli and let
$a\uparrow b(\beta-1)$ along a countable sequence to obtain
\eqref{eq:laplace-gap-lower}.
\end{proof}

\section{Further questions}
\label{sec:questions}

The following questions concern information not supplied by the preceding
proofs.
Throughout the first three questions the offspring variance is finite
and positive and $\mu$ is nonconstant.

\subsection{The constant in the one-time probability}
Does
\[
 n\Q(A_n)\longrightarrow\gamma(p,\mu)\in(0,\infty)
\]
hold? Since $a_n=q_n\E[K_n\mid Z_n>0]$, this is equivalent to the
existence of the limit of the conditional mean multiplicity.
For atomless $\mu$ the answer is already $\gamma=2/\sigma^2$.
A first discrete case is $p_0=p_2=1/2$ with
$\mu=(\delta_{-1}+\delta_1)/2$.
The uniform conditional second moment proved here gives tightness and
uniform integrability of the conditional multiplicities. Thus a unique
conditional distributional limit for $K_n$ would imply the desired
constant; identifying that limit remains a separate task.

\subsection{Fluctuations at the logarithmic scale}
For a specified nonconstant law, does $D_N/\log N$ converge in
probability, almost surely, or only in distribution? Even when $\mu$
is atomless, the identity $\EQ[D_N]\sim(2/\sigma^2)\log N$ does not
answer this question. Natural initial cases are the preceding binary
nearest-neighbor model and binary branching with standard Gaussian
increments. The constant-displacement example in
Section~\ref{sec:fluctuations} concerns a different, selected event
and does not settle these cases.

A concrete route is to sharpen the multipoint bound
\eqref{eq:generic-multipoint} to asymptotics in regimes contributing to
factorial moments. Two-time upper bounds alone lose the effect of old
bushes which remain alive and compete with the spine subtree. Quantifying
that competition is the missing spatial input.

\subsection{Arrival times and separation from the minimum}
The arrival-time conclusion is equivalently
$\log\tau_k=k^{1+o(1)}$. Can one determine fluctuations of
$\log\tau_k$, or the sizes of the gaps between consecutive extremal
generations? Another observable is
\[
 G_n=V(w_n)-M_n\geq0.
\]
Under centered finite-variance spatial increments and finite offspring
variance, one can ask for scaling limits of $G_n$ and for counts of the
near-minimum events $\{G_n\leq h_n\}$. Fixed-time extremal-displacement
limit theorems do not by themselves determine these temporal counts.
In the infinite-variance family of Theorem~\ref{thm:finite-main},
Theorem~\ref{thm:quantitative-gap} gives eventual lower bounds for
Gaussian and Laplace increments; their sharpness remains undetermined.

\subsection{Other infinite-variance offspring laws}
For the explicit stable family of Theorem~\ref{thm:finite-main},
nonconstant displacements give finitely many weak-minimum events,
including in the lattice case. For more general critical offspring
laws, the selected event has the complete survival-series criterion
of Theorem~\ref{thm:criterion}; when $\mu$ is atomless, it is also the
criterion for $A_n$. With ties, the generic counting argument reduces
recurrence exactly to divergence of $\sum_n\E[K_n]$.
A remaining question is whether $\sum_nq_n<\infty$ can coexist with
weak-minimum recurrence for a nonconstant displacement law outside the
stable family considered here. A further question is to replace the
exact stable generating function by general regularly varying
infinite-variance offspring tails while retaining a useful multiplicity
comparison.

\section*{Acknowledgement}
The author thanks Antal J\'arai for proposing the problem that motivated
this work.


\clearpage
\begin{thebibliography}{9}

\bibitem{EP17}
J. Engl\"ander and Y. Peres,
\emph{Survival asymptotics for branching random walks in IID environments},
Electron. Commun. Probab. \textbf{22} (2017), paper no.~29, 1--12.
\href{https://doi.org/10.1214/17-ECP60}{doi:10.1214/17-ECP60}.


\bibitem{JN14}
A. A. J\'arai and A. Nachmias,
\emph{Electrical resistance of the low dimensional critical branching
random walk},
Comm. Math. Phys. \textbf{331} (2014), 67--109.
\href{https://doi.org/10.1007/s00220-014-2085-y}{doi:10.1007/s00220-014-2085-y}.

\bibitem{K86}
H. Kesten,
\emph{Subdiffusive behavior of random walk on a random cluster},
Ann. Inst. H. Poincar\'e Probab. Statist. \textbf{22} (1986), no.~4,
425--487.
\href{https://www.numdam.org/item/AIHPB_1986__22_4_425_0/}{Numdam full text}.

\bibitem{K95}
H. Kesten,
\emph{Branching random walk with a critical branching part},
J. Theoret. Probab. \textbf{8} (1995), no.~4, 921--962.
\href{https://doi.org/10.1007/BF02410118}{doi:10.1007/BF02410118}.

\bibitem{KNS66}
H. Kesten, P. Ney and F. Spitzer,
\emph{The Galton--Watson process with mean one and finite variance},
Theory Probab. Appl. \textbf{11} (1966), no.~4, 513--540.
\href{https://doi.org/10.1137/1111059}{doi:10.1137/1111059}.

\bibitem{KS64}
S. Kochen and C. Stone,
\emph{A note on the Borel--Cantelli lemma},
Illinois J. Math. \textbf{8} (1964), no.~2, 248--251.
\href{https://doi.org/10.1215/ijm/1256059668}{doi:10.1215/ijm/1256059668}.


\bibitem{LS15}
S. P. Lalley and Y. Shao,
\emph{On the maximal displacement of critical branching random walk},
Probab. Theory Related Fields \textbf{162} (2015), nos.~1--2, 71--96.
\href{https://doi.org/10.1007/s00440-014-0566-8}{doi:10.1007/s00440-014-0566-8}.
\href{https://arxiv.org/abs/1212.2933}{arXiv:1212.2933}.

\bibitem{LPP95}
R. Lyons, R. Pemantle and Y. Peres,
\emph{Conceptual proofs of $L\log L$ criteria for mean behavior of
branching processes},
Ann. Probab. \textbf{23} (1995), no.~3, 1125--1138.
\href{https://doi.org/10.1214/aop/1176988176}{doi:10.1214/aop/1176988176}.

\bibitem{Slack68}
R. S. Slack,
\emph{A branching process with mean one and possibly infinite variance},
Z. Wahrscheinlichkeitstheorie verw. Gebiete \textbf{9} (1968), 139--145.
\href{https://doi.org/10.1007/BF01851004}{doi:10.1007/BF01851004}.

\end{thebibliography}
\end{document}